\documentclass{article}
\usepackage{amssymb,amsmath}

\newtheorem{theorem}{Theorem}[section]

\newtheorem{proposition}[theorem]{Proposition} 

\numberwithin{equation}{section}

\newcommand{\ph}{\ensuremath{\varphi}}
\newcommand{\ex}[1]{\ensuremath{\mathrm{e}^{{#1}}}}
\newcommand{\exi}[1]{\ensuremath{\ex{\eul {#1}}}}
\newcommand{\exni}[1]{\ensuremath{\ex{-\eul {#1}}}}
\newcommand{\dif}[1]{\ensuremath{\:\mathrm{d}{#1}}}
\newcommand{\reel}{\ensuremath{\mathbb{R}}}
\newcommand{\cpl}{\ensuremath{\mathbb{C}}}
\newcommand{\Kfrac}[2]{\ensuremath{\displaystyle\frac{{#1}}{{#2}}}}
\newcommand{\nm}[1]{\ensuremath{\|{#1}\|}}
\newcommand{\eul}{\ensuremath{\mathrm{i}}}
\newcommand{\ens}[2]{\ensuremath{\{{#1}\;;\;{#2}\}}}
\newcommand{\eps}{\ensuremath{\varepsilon}}
\newcommand{\KInv}[1]{\ensuremath{\displaystyle\Inv{{#1}}}}
\newcommand{\Inv}[1]{\ensuremath{\frac{1}{{#1}}}}
\newcommand{\fml}[3]{\ensuremath{\lbrace {#1} \rbrace_{{#2}}^{{#3}}}}
\newcommand{\rel}{\ensuremath{\mathbb{Z}}}
\newcommand{\sgl}[1]{\ensuremath{\{{#1}\} }}
\newcommand{\dem}{\ensuremath{\frac{1}{2}}}
\newcommand{\Hil}{\mathcal{H}}
\newcommand{\conj}[1]{\overline{{#1}}}
\newcommand{\ps}[2]{\langle{#1},{#2}\rangle}
\newcommand{\Nev}{\mathcal{N}}
\newcommand{\Poisson}{\mathcal{P}}
\newcommand{\dpar}[1]{\frac{\partial}{\partial{#1}}}
\newcommand{\Dpar}[2]{\frac{\partial{#2}}{\partial{#1}}}
\newcommand{\dparh}[2]{\frac{\partial^{#2}}{\partial{#1}^{#2}}}
\newcommand{\Dparh}[3]{\frac{\partial^{#2}{#3}}{\partial{#1}^{#2}}}
\newcommand{\Alpha}{\mathcal{A}} 
\newcommand{\dist}[2]{\ensuremath{\mathrm{dist}({#1},{#2})}}

\newcommand{\kfrac}[2]{\ensuremath{{{#1}}/{{#2}}}}

\newcommand{\ff}{\textstyle\frac{5}{4}\displaystyle}

\begin{document}

\title{de Branges spaces with bi-Lipschitz phase for large distances}

\author{Philippe Poulin \thanks{United Arab Emirates University, Al Ain, U.A.E., philippepoulin@uaeu.ac.ae} 
        \and Simon Cowell \thanks{United Arab Emirates University, Al Ain, U.A.E., scowell@uaeu.ac.ae}}

\maketitle

\begin{abstract}
 We characterize an arbitrary de Branges space with bi-Lipschitz phase for large distances as a subspace of a weighted
  Paley--Wiener space, consisting of the elements square-integrable against an explicit extra-weight on the real line.
\end{abstract}

\section{Introduction}



The inverse Fourier transform, defined by
 $$\mathcal{F}^{-1}[\ph](z)=\Inv{\sqrt{2\pi}}\int_{-\pi}^\pi\exi{zt}\ph(t)\dif{t}, \ \ \ z\in\cpl$$
maps $L^2[-\pi,\pi]$ to a space of entire functions known as the classical \emph{Paley--Wiener space} and denoted by $L^2_\pi$.
By the Paley--Wiener theorem, $L^2_\pi$ consists of all entire functions of exponential type at most $\pi$ that are square-integrable on $\reel$.
In other words,
  $$L^2_\pi=\ens{f \mbox{ entire }}{\nm{f}_2<\infty, \ |f(z)|\leq C_\eps\ex{(\pi+\eps)|z|}},$$
where the norm, $\nm{\ }_2$, is induced by the usual Hermitian product $\ps{f}{g}_2=\displaystyle\int_{-\infty}^\infty f(t)\conj{g(t)}\dif{t}$. 

The Hilbert space $L^2_\pi$ admits a \emph{reproducing kernel} at each $\zeta\in\cpl$, that is, a function $k_\zeta\in L^2_\pi$ satisfying
 $\ps{f}{k_\zeta}_2=f(\zeta)$. 
Furthermore, it admits an \emph{orthonormal basis} of reproducingkernels, namely $\fml{k_n}{n\in\rel}{}$.
Explicitly,
 $$k_\zeta(z)=\frac{\sin(\pi(z-\bar{\zeta}))}{\pi(z-\bar{\zeta})},$$
yielding the following expansion \cite[p.150]{L}: for $f\in L^2_\pi$,
 $$f(z)=\sum_{n=-\infty}^\infty d_n\frac{\sin(\pi(z-n))}{\pi(z-n)}.$$


In \cite{dB}, de Branges has identified a large class of spaces of entire functions which also admit orthogonal bases of reproducing kernels \cite[p.55]{dB} and can be seen as a generalization of $L^2_\pi$.
A Hilbert space $\Hil$ of entire functions is a \emph{de Branges space} if it satisfies the following properties \cite[p.57]{dB}:
\emph{%
 \begin{enumerate}
  \item
   The linear functional $\Hil\to\cpl,\ f\mapsto f(z_0)$ is bounded for all $z_0\in\cpl$;
  \item
   If $f(z)\in\Hil$, then $f^*(z)$ also belongs to $\Hil$ and has the same norm as $f(z)$, where $f^*(z)=\conj{f(\bar{z})}$;
  \item
   If $f(z)\in\Hil$ and $f(z_0)=0$, then $f(z)\Kfrac{z-\conj{z_0}}{z-z_0}$ also belongs to $\Hil$ and has the same norm as $f(z)$.
 \end{enumerate}
}

\noindent By the Riesz lemma, the first property ensures that $\Hil$ admits a reproducing kernel at each $\zeta\in\cpl$.

Concrete examples of de Branges spaces may be produced as follows.
A function $f(z)$ analytic in $\cpl^+$ is in the \emph{Nevanlinna class} if it is the ratio of two analytic, bounded functions in $\cpl^+$.
The \emph{mean type} of such a function is then given by
 $$\limsup_{y\to\infty}\frac{\log|f(\eul y)|}{y}.$$
For $h\in\reel$, let us denote by $\Nev_h^+$ the class of functions in the Nevanlinna class whose mean type does not exceed $h$.  
Given a \emph{Hermite--Biehler function}, that is, an entire function such that $|E(\bar{z})|<|E(z)|$ for $\Im z>0$, 
 $$\Hil(E)=\ens{f\mbox{ entire}}{\nm{f/E}_2<\infty\mbox{ and }f/E,\ f^*/E\in\Nev_0^+}$$
is a de Branges space, equipped with the norm $\nm{f}=\nm{f/E}_2$. 
Indeed, de Branges proved that \emph{each de Branges space is isometrically equal to a space of the form $\Hil(E)$} \cite[p.57]{dB}.
Notice that $L^2_\pi$ itself is a de Branges space, obtained from the Hermite--Biehler function $E(z)=\exni{\pi z}$.\\


In their work about sampling and interpolation, Lyubarskii and Seip studied a large class of de Branges spaces that includes
 $L^2_\pi$, defined as follows \cite{LS}. 
Given a de Branges space $\Hil$, let $M(z)=\nm{k_z}$.
$\Hil$ is a \emph{weighted Paley--Wiener space} if the restriction of $M$ to the real axis satisfies
 \begin{enumerate}
  \item
   $M(x)>0$ for all $x\in\reel$;
  \item
   $\nm{f}\simeq \nm{f/M}_2$ for all $f\in\Hil$.\,\footnote{%
In the sequel we say that two positive functions $f$ and $g$ are \emph{comparable} and write $f\simeq g$ if there exist 
constants $A$, $B>0$ such that $Af\leq g\leq Bf$.}
 \end{enumerate}

Concrete examples of weighted Paley--Wiener spaces are produced as follows.
Let $m(t)\simeq 1$ be a measurable function on $\reel$.
The \emph{potential} of the measure $m(t)\dif{t}$ is defined as
 $$\omega_m(z)=\int_{-\infty}^\infty\log^*\left| 1-\frac{z}{t}\right|m(t)\dif{t},$$
where 
 \begin{equation*}
  \log^*|1-z/t|=
   \left\{
    \begin{array}{ll}
     \log|1-z/t|+(\Re z)/t&\mbox{ if } |t|> 1\\
     \log|1-z/t|&\mbox{ otherwise.}
    \end{array}
   \right.
 \end{equation*}

\noindent Then, 
 $$PW(m)=\ens{f\mbox{ entire}}{\nm{f\ex{-\omega_m}}_2<\infty\mbox{ and }|f(z)|\ex{-\omega_m(z)}\leq C_\eps\ex{\eps|\Im z|}}$$
is a weighted Paley--Wiener space.
For $g(z)$ real-entire, $\ex{g}PW(m)$ is also a weighted Paley--Wiener space.
Indeed, Lyubarskii and Seip proved that \emph{any weighted Paley--Wiener space is equal, with equivalence of norms, to a space of the form $\ex{g}PW(m)$}.
As an example, $PW(1)$, the simplest weighted Paley--Wiener space, is equal to $L^2_\pi$.

\paragraph{Motivation and main result}

Let $\Hil$ be a de Branges space and let $k_z\in\Hil$ be the reproducing kernel at $z\in\cpl$.
The sequence of complex number $\sgl{z_j}$ is \emph{interpolating} for $\Hil$ if there exists an $f\in\Hil$ satisfying 
 $f(z_j)=a_j$ for any choice of interpolation data $\sgl{a_j/\nm{k_{z_j}}}\in\ell^2(\cpl)$ \cite[p.21]{S}.
It is \emph{complete interpolating} if in addition $f$ is unique.
The sequence $\sgl{z_j}$ is \emph{sampling} if 
 $\nm{f}^2\simeq \sum |f(z_j)|^2/\nm{k_{z_j}}^2.$
The notions of \emph{sampling} and \emph{interpolating} sequences are often presented as dual.
It is well known \cite[p.3]{MNO} that $\sgl{z_j}$ is complete interpolating if and only if it is interpolating and sampling.

Classes of de Branges spaces are often studied in order to understand their sampling and interpolating sequences \cite{LS,S,MNO,BMS,B}.
Results are obtained for weighted Paley--Wiener spaces \cite{LS}.
Notice that the Hermite--Biehler function of such a space may be chosen with zeroes \emph{equidistributed} on the same line, in the sense
 that the distance between two consecutive zeroes is comparable with $1$. 
At the other extreme, results are also obtained for de Branges spaces coming from a Hermite--Biehler function with sparse zeroes \cite{BMS}.

The question is asked \cite[p.5]{MNO} to identify a class of de Branges spaces broader than the weighted Paley--Wiener spaces,
 but still opposite to the sparse case, in which some structural results may be obtained.
Indeed, a class of de Branges spaces larger than the weighted Paley--Wiener spaces is already studied in \cite{LS}.
Its definition and the results implicitly present in \cite{LS} are given in Section~\ref{Sec:Examples}; see \cite{P} for a complete exposition.

The present paper proposes to study an even larger class of de Branges spaces, defined as follows.
Let $\Hil(E)$ be a de Branges space whose Hermite--Biehler function $E$ does not have any real zeroes. 
Then $E(x)$ admits a polar decomposition $|E(x)|\exni{\ph(x)}$ on the real axis, where $\ph(x)$, the so-called \emph{phase}, is real-analytic and well-defined up to the addition of $2k\pi$.
It is well known that $\ph(x)$ is also increasing \cite[p.54]{dB}.
We are interested in the case where $\ph$ is \emph{bi-Lipschitz for large distances}: there exist positive constants $N$, $C_1$, and $C_2$ such that
 $$C_1(x_2-x_1)\leq\ph(x_2)-\ph(x_1)\leq C_2(x_2-x_1) \mbox{\ \ whenever }x_2-x_1> N,$$
where $C_1$, $C_2$, and $N$ are independent of $x_1$ and $x_2$.

Let us state our main result.
Following Lyubarskii and Seip \cite{LS}, we shall say that $f(z)$ is of \emph{$\omega_m$-type} if for each $\eps>0$, there exists a $C_\eps>0$ such that
 \begin{equation*}
 |f(z)|\leq C_\eps\ex{\eps|z|}\ex{\omega_m(z)}
 \end{equation*}
in the complex plane.

\begin{theorem}
 Let $\Hil=\Hil(E)$ be a de Branges space, with $E(x)=|E(x)|\exni{\ph(x)}$ on the real line.
 Assume $\ph$ is bi-Lipschitz for distances larger than $N\geq 0$. 
 Then, there exists a measurable $m\simeq 1$ and a real-entire $g(z)$ such that
  $$\Hil\subseteq \ex{g}PW(m).$$
 Namely, 
  $$\ex{-g}\Hil=\ens{f\mbox{ \rm{entire}}}{\nm{f(x)\ex{-\omega_m(x)}\ex{\theta(x)}}_2<\infty,\ |f(z)|\ex{-\omega_m(z)}\leq C_\eps \ex{\eps|\Im z|}},$$
 where $\theta(x)=\KInv{\pi}\int_{|t-x|<N}\frac{\ph(t)-\ph(x)}{t-x}\dif{t}$. 
\end{theorem}

In their study of weighted Paley--Wiener spaces, Lyubarskii and Seip developed a \emph{multiplier lemma} in the line of Beurling \cite{B}.
It produces, for $\tau>\sup m$, a Hermite--Biehler function $E_{\tau-m}$ whose zeroes are equidistributed along the axis $\Im z=-1$.
Let $\Sigma_{\tau-m}$ be the set of zeroes of $E_{\tau -m}$.
The multiplication by $E_{\tau-m}$ is then a bijection with equivalence of norms from $PW(m)$ to 
 $$L^2_{\pi\tau}[\Sigma_{\tau-m}]=\ens{f\mbox{ entire}}{\nm{f(x)}_2<\infty,\ |f(z)|<C_\eps\ex{(\pi\tau+\eps)|\Im z|},\ f(\Sigma_{\tau-m})=0}.$$
In this way, Lyubarskii and Seip tranferred results about sampling and interpolation from $L^2_{\pi\tau}[\Sigma_{\tau-m}]$ to $PW(m)$.
Our theorem implies that a similar transfer may be done from
 $$\ens{f\mbox{ entire}}{\nm{f(x)\ex{\theta(x)}}_2<\infty,\ |f(z)|<C_\eps\ex{(\pi\tau+\eps)|\Im z|},\ f(\Sigma_{\tau-m})=0}$$
 to $\ex{-g}\Hil$.

To prove our theorem we shall follow part of Lyubarskii and Seip's study of weighted Paley--Wiener spaces with two modifications. 
Firstly, we shall generalize the definition of $\omega_\gamma$ so that $\gamma(x)$ may be any positive, continuous function whose
 antiderivative $\ph(x)$ satisfies $|\ph(x)-\ph(0)|\simeq |x|$.
Efforts are made to show that $\omega_\gamma$ satisfies the expected properties of the potential of a measure in $\cpl$ (see \cite{R}).
Secondly, we shall generalize Lyubarskii and Seip's multiplier lemma so the zeroes of the resulting Hermite--Biehler function may be multiple.

\paragraph{Notation and terminology} 

For positive functions $f$ and $g$, we write $f\lesssim g$ if there exists a constant $C$ such that $f\leq Cg$ pointwise.
As already mentioned, we write $f\simeq g$ and say that \emph{$f$ is comparable with $g$} if $f\lesssim g$ and $g\lesssim f$.
We also say that a set of points is \emph{equidistributed} if the distance between a point and its closest neighbor is
 comparable with $1$. 

For real-valued functions $f$ and $g$, we write $f\gg g$ if $f>g+\eps$ for a certain $\eps>0$.

For a complex-valued function $f$, we define $f^*(z)=\conj{f(\bar{z})}$.
The $L^2$-norm of $f(z)$ is defined as $\nm{f}_2^{\,2}=\int_{-\infty}^\infty |f(x)|^2\dif{x}$, provided that this last integral converges.
Finally, we say that $f(z)$ is \emph{real-entire} if it is an entire function real on the real line.

\section{Proofs}

\subsection{Potentials considered}

Let $\ph(x)$ be a continuously differentiable, increasing function satisfying $|\ph(x)|\simeq |x|$ for $|x|$ large. 
Our aim is to verify that
 $$\omega_{\ph'}(z)=\int_{-\infty}^\infty\log^*\left|1-\frac{z}{t}\right|\ph'(t)\dif{t}$$
satisfies the expected properties of a potential.

Let us show that the above integral is absolutely convergent.
The condition on $\ph$ ensures that for $R>|z|+1$, 
 $$\int_{|t|>R}\dpar{t}\left(\log^*\left|1-\frac{z}{t}\right|\right)\ph(t)\dif{t}=\int_{|t|>R}\Re\left(\frac{z^2}{t^2(t-z)}\right)\ph(t)\dif{t}$$
is well-defined. 
The following integral is thus also well-defined,
 $$\int_{|t|>R}\log^*\left|1-\frac{z}{t}\right|\ph'(t)\dif{t}= A_R(z)-\int_{|t|>R}\left(\dpar{t}\log^*\left|1-\frac{z}{t}\right|\right)\ph(t)\dif{t},$$
where $A_R(z)=\log^*|1+z/R|\ph(-R)-\log^*\left|1-z/R\right|\ph(R)$.
Since its integrand changes sign finitely many times, it is absolutely convergent.
The result follows.

The continuity of $\omega_{\ph'}$ also follows from the previous relation, by applying the dominated convergence theorem to the integral on its right-hand side.

Let us show that for $y\neq 0$, $z=x+\eul y$, $\dpar{y}\omega_{\ph'}$ may be calculated by interchanging the derivative and the integral.
The dominated convergence theorem implies
 $$\dpar{y}\int_{|t|>R}\log^*\left|1-\frac{z}{t}\right|\ph'(t)\dif{t} = \dpar{y}A_R(z)-\int_{|t|>R}\dpar{t}\left(\dpar{y}\log^*\left|1-\frac{z}{t}\right|\right)\ph(t)\dif{t}.$$
The last integral in the above relation may be evaluated by parts, yielding %
 $$\dpar{y}\int_{|t|>R}\log^*\left|1-\frac{z}{t}\right|\ph'(t)\dif{t}=\int_{|t|>R}\dpar{y}\log^*\left|1-\frac{z}{t}\right|\ph'(t)\dif{t}.$$
The dominated convergence theorem may also be used for the interval of integration $[-R,R]$, yielding in total
 $$\dpar{y}\omega_{\ph'}(z)=\int_{-\infty}^\infty\dpar{y}\log^*\left|1-\frac{z}{t}\right|\ph'(t)\dif{t}=\int_{-\infty}^\infty\frac{y}{(x-t)^2+y^2}\ph'(t)\dif{t}=\pi \Poisson_{\ph'}(z),$$
where $\Poisson$ denotes the Poisson transform.

Finally, let us point out that the distributional Laplacian of $\omega_{\ph'}$ is given by 
 $$\Delta\omega_{\ph'}(x+\eul y)=2\pi\ph'(x)\dif{x}\dif{\delta_0}(y),$$ 
where $\delta_0$ denotes the Dirac measure at $0$. 
The proof is exactly the same as when $\ph'(x)\simeq 1$, see \cite[p.6]{P}, which follows closely \cite[p.74]{R}.
 
\subsection{Multiplier lemma} 

We now develop a version of Lyubarskii and Seip's multiplier lemma dealing with a positive, continuous $\gamma(t)$ whose antiderivative is 
 bi-Lipschitz for large distances. 
We shall obtain an equivalence $\ex{\omega_\gamma(z)}\simeq |F_\gamma(z)|$ for a real-entire function $F_\gamma(z)$ whose zeroes are equidistributed on the real line, but now have multiplicity. 
It will not be possible to shift them to $\cpl^-$ without breaking the equivalence, which thus holds on an upper half-plane $\Im z>\eps>0$ only.  

Let $\ph(x)=\int_0^x\gamma(t)\dif{t}$.
By hypothesis, there exist $C_1$, $C_2$, and $N$ such that 
 \begin{equation}\label{BiLipschitz}
  C_1(x_2-x_1)\leq \ph(x_2)-\ph(x_1)\leq C_2(x_2-x_1) \mbox{\ \ whenever } x_2-x_1\geq N.
 \end{equation}
Consider a sequence $\Alpha=\fml{\alpha_k}{k\in\rel}{}$ of natural numbers satisfying $2C_2N<\alpha_k<B$ for a certain bound $B$. 
Let $\cdots<x_{-1}<x_0<x_1<\cdots$ be the partition of $\reel$ defined by
 $x_0=0$ and $\int_{x_k}^{x_{k+1}}\gamma(t)\dif{t}=\alpha_k$.

We claim that $x_{k+1}-x_k\simeq 1$.
Indeed, for $x_k\leq x\leq x_k+2N$
 $$\ph(x)-\ph(x_k)\leq \ph(x_k+2N)-\ph(x_k)\leq 2C_2N<\alpha_k, $$
which yields $x_{k+1}>x_k+2N$.
Similarly, for $x-x_k\geq B/C_1$
 $$\ph(x)-\ph(x_k)\geq B>\alpha_k,$$
which yields $x_{k+1}<x_k+(B/C_1)$.
Therefore, $2N<x_{k+1}-x_k<B/C_1$, as desired.

Let $\xi_k=\KInv{\alpha_k}\int_{x_k}^{x_{k+1}}t\gamma(t)\dif{t}$.
We claim that $\xi_k$, which lies in $(x_k,x_{k+1})$, is bounded away from $x_k$, and hence $\xi_{k+1}-\xi_k\simeq 1$.
In fact,
 \begin{equation*}
  \xi_k\geq\Inv{\alpha_k}\int_{x_k}^{x_k+N}\hspace{-0.8cm} x_k \gamma(t)\dif{t} + \Inv{\alpha_k}\int_{x_k+N}^{x_{k+1}}\hspace{-0.6cm} (x_k+N)\gamma(t)\dif{t}
   = x_k+\frac{N}{\alpha_k}(\ph(x_{k+1})-\ph(x_k+N)).
 \end{equation*}
Moreover, $\ph(x_{k+1})-\ph(x_k+N)\geq C_1(x_{k+1}-x_k-N)>C_1 N$.
Therefore, $\xi_k-x_k\geq C_1N^2/\alpha_k\geq C_1N^2/B$, as desired.

Consider the auxiliary measure
 $$\dif{\nu}(t)=\gamma(t)\dif{t}-\sum_{k=-\infty}^{\infty}\alpha_k\dif{\delta_{\xi_k}}(t),$$
where $\dif{\delta_\xi}$ denotes the Dirac measure at $\xi\in\reel$.
Let $f_\nu(x)=\int_0^x\dif{\nu}(t)$ (inclusive of the endpoints) and $g_\nu(x)=\int_0^x f_\nu(t)\dif{t}$.
For $a<b$ and $z\notin\reel$, two integrations by parts give
 \begin{equation*}
  \int_{[a,b]}\log|1-z/t|\dif{\nu}(t) = \int_a^b g_\nu(t)\dparh{t}{2}\log|1-z/t|\dif{t} + R(a,b),
 \end{equation*}
where $R(a,b)=(f_\nu(t)\log|1-z/t|-g_\nu(t)\dpar{t}\log|1-z/t|)\ |_a^b$.

Since $\xi_k\in(x_{k},x_{k+1})$, clearly $f_\nu(x_k)=0$ for all $k$.
Since in addition $\Alpha$ is bounded, $f_\nu$ is a bounded function on $\reel$.
Moreover, $g_\nu(x_k)=0$ for all $k\in\rel$ by choice of $\xi_k$.
In addition, $g_\nu$ is bounded on $\reel$.

We deduce on the one hand that $R(a,b)\to 0$ when $a\to-\infty$ and $b\to\infty$.
On the other hand, using any appropriate branch of the logarithm,
 \begin{eqnarray*} \int_a^b g_\nu(t)\dparh{t}{2}\log|1-z/t|\dif{t}
  &=& \int_a^b g_\nu(t)\,\Re\dparh{t}{2}\log(1-z/t)\dif{t}\\
  &=& -\int_a^b g_\nu(t)\Re \Inv{(t-z)^2}\dif{t}+\int_a^b\frac{g_\nu(t)}{t^2}\dif{t}.
 \end{eqnarray*}
Since $g_\nu$ is bounded, and since $g_\nu(t)$ is nonnegative and $\simeq t^2$ in a neighborhood of $x_0=0$, we conclude that
 $\int_{-\infty}^\infty\log|1-z/t|\dif{\nu}(t)$
is well-defined and satisfies
 \begin{equation*}
  \left|\int_{-\infty}^\infty\log|1-z/t|\dif{\nu}(t)\right|\lesssim 1
 \end{equation*}
when $\Im z\gg 0$ (i.e., when $\Im z>\eps>0$ for a certain $\eps>0$).

An integration by parts also gives
 $$\int_{-\infty}^\infty\frac{\chi_{\reel\setminus[-1,1]}(t)}{t}\dif{\nu(t)}=-f_\nu(-1)-f_\nu(1)+\int_{|t|>1}\frac{f_\nu(t)}{t^2}\dif{t}.$$
Since $f_\nu$ is bounded, the above expression is just a real constant $C$.
In total,
 \begin{equation*}
  \left|\int_{-\infty}^\infty\log^*|1-z/t|\dif{\nu}(t)-Cx\right|\lesssim 1
 \end{equation*}
when $\Im z\gg 0$.
Letting $\alpha=C-\sum_{|\xi_k|\leq 1}\alpha_k/\xi_k$, it follows that
 $$|\,\omega_\gamma(z)-\alpha x-\textstyle\sum_k\,\alpha_k(\log|1-z/\xi_k|+x/\xi_k)\,|\lesssim 1.$$
In other words, for $F_\gamma(z)=\ex{\alpha z}\prod_k(1-z/\xi_k)^{\alpha_k}\ex{z\,\alpha_k/\xi_k}$,
 $$|F_\gamma(z)|\simeq\ex{\omega_\gamma(z)}$$
when $\Im z\gg 0$.\\

We have proven the following \emph{multiplier lemma}:

\begin{proposition}\label{Prop:ML:ML}
 Let $\gamma(x)$ be a positive, continuous function whose antiderivative $\ph(x)$ satisfies the condition (\ref{BiLipschitz}).
 Let $\Alpha=\fml{\alpha_k}{k\in\rel}{}$ be a bounded sequence of positive integers larger than $2C_2N$.
 Then, $F_\gamma(z)=\ex{\alpha z}\prod_k(1-z/\xi_k)^{\alpha_k}\ex{z\,\alpha_k/\xi_k}$ is a real-entire function satisfying
  $$|F_\gamma(z)|\simeq\ex{\omega_\gamma(z)}\ \mbox{ when }\Im z\gg 0.$$
 The consecutive zeroes $\sgl{\xi_k}$ of $F_\gamma$ are real, have respective multiplicity $\alpha_k$, and satisfy $\xi_{k+1}-\xi_k\simeq 1$.
\end{proposition}

\subsection{Proof of the theorem}

Let $\Hil=\Hil(E)$ be a de Branges space whose phase $\ph$ satisfies the equation (\ref{BiLipschitz}). 
Since $E$ is a Hermite--Biehler function, $\log|E(x+\eul|y|)|$ is subharmonic. 
The computation of its distributional Laplacian is the same as that for weighted Paley--Wiener spaces.
Indeed, for $H(x,y)=\log|E(x+\eul y)|$, $\dparh{x}{2}(H(x,|y|))=(\dparh{x}{2}H)(x,|y|)$, while
 $$\dparh{y}{2}(H(x,|y|))=\Dparh{y}{2}{H}(x,|y|)+2\Dpar{y}{H}(x,|y|)\dif{\delta_0(y)}\dif{x}$$
in the sense of distributions, where $\delta_0$ denotes the Dirac measure at $0$. 
The harmonicity of $H$ in the neighborhood of the closed upper half-plane implies
 $$\Delta(H(x,|y|))=(\Delta H)(x,|y|)+2\Dpar{y}{H}(x,|y|)\dif{\delta_0(y)}\dif{x}=2\Dpar{y}{H}(x,0)\dif{x}\dif{\delta_0(y)}.$$
Since $\log E(x)=H(x,0)-\eul\ph(x)$ has an analytic extension in the neighborhood of $\reel$, the Cauchy--Riemann
 equations yield $\Delta(H(x,|y|))=2\ph'(x)\dif{x}\dif{\delta_0(y)}$.
Consequently,
 $$\log|E(x+\eul|y|)|=h(z)+\omega_{\ph'/\pi}(z),$$
where $z=x+\eul y$ and $h$ is harmonic.
Since $h(\bar{z})=h(z)$, $h$ is indeed the real part of a real-entire function $g(z)$.
It follows that for $\gamma=\ph'/\pi$
 $$\ex{-g}\Hil=\ens{f\mbox{ entire}}{\nm{f\ex{-\omega_\gamma}}_2<\infty,\ f^\#\ex{-\omega_\gamma-\eul\tilde{\omega}_\gamma} \in \Nev_0^+},$$
where $f^\#$ varies in $\sgl{f,f^*}$.

Observe that $\omega_\gamma(z+\eul)>\omega_\gamma(z)$ on the upper half-plane. 
In fact, $\ex{\omega_\gamma(z)-\omega_\gamma(z+\eul)}$ and $\ex{\omega_\gamma(z+\eul)-\omega_\gamma(z)}$ are moduli of 
 functions in $\Nev_0^+$, since the mean value theorem and the dominated convergence theorem imply
  $$0 \leq \limsup_{y\to\infty} \frac{\omega_\gamma(\eul (y+1))-\omega_\gamma(\eul y)}{y} 
	    \leq \limsup_{y\to\infty}\frac{y+1}{y}\int_{-\infty}^\infty\Inv{t^2+y^2}\gamma(t)\dif{t}=0.$$
Therefore, $\omega_\gamma(z)$ may be replaced with $\omega_\gamma(z+\eul)$ in the condition
 $f^\#\ex{-\omega_\gamma-\eul\tilde{\omega}_\gamma}\in\Nev_0^+$. 
Moreover, the multiplier lemma ensures that $\ex{\omega_\gamma(z)}\simeq F_\gamma(z)$ on $\Im z\gg 0$ for an entire function $F_\gamma(z)$ with 
 equidistributed (but multiple) zeroes on the real axis.
Obviously, there exists an $m(x)\simeq 1$ generating the same zeroes with the same multiplicities, so
 $\ex{\omega_m(z)}\simeq F_m(z) = F_\gamma(z)\simeq \ex{\omega_\gamma(z)}$ when $\Im z\gg 0$.%
\footnote{%
In practice, any choice of a concrete $m$ such that $\ex{\omega_m(z)}\simeq \ex{\omega_\gamma(z)}$ when $\Im z\gg 0$ works.
} 
In particular, $\ex{\omega_\gamma(z+\eul)}\simeq \ex{\omega_m(z+\eul)} \simeq \ex{\omega_m(z)}$ for $\Im z\geq 0$, where the last equivalence is 
 justified by the mean value theorem.
In total, the condition $f^\#\ex{-\omega_\gamma-\eul\tilde{\omega}_\gamma}\in\Nev_0^+$ may thus be replaced with
 $f^\#\ex{-\omega_m-\eul\tilde{\omega}_m}\in\Nev_0^+$ in the characterization of $\ex{-g}\Hil$. 

As for weighted Paley--Wiener spaces, we may replace this last condition with an omega type condition.
In fact, let $\tau>\sup m$ and $E_{\tau-m}(z)=F_{\tau-m}(z+\eul)$, where $F_{\tau-m}$ is given by the multiplier lemma (with zeroes of
 multiplicity $1$). 
Then, $|E_{\tau-m}(z)|\simeq\ex{\omega_{\tau-m}(z)}$ for $\Im z\geq 0$, while $\omega_\tau(x+\eul y)=\pi\tau|y|$.
The condition $f^\#\ex{-\omega_m-\eul\tilde{\omega}_m}\in\Nev_0^+$ is thus equivalent to $f^\# E_{\tau-m}\in\Nev_{\pi\tau}^+$.
Observe that $E_{\tau-m}$ is a Hermite--Biehler function satisfying $E_{\tau-m}(\bar{z})=E_{\tau-m}(z-2\eul)$, and hence the mean value theorem implies
 $$\left|\frac{E_{\tau-m}^*(\eul y)}{E_{\tau-m}(\eul y)}\right|\simeq\ex{\omega_{\tau-m}(\eul (y-2))-\omega_{\tau-m}(\eul y)}\simeq 1.$$
In particular, $E^*_{\tau-m}/E_{\tau-m}\in\Nev_0^+$.
The condition $f^\# E_{\tau-m}\in\Nev_{\pi\tau}^+$ is thus equivalent to $f^\#E_{\tau-m}\in\Nev_{\pi\tau}^+$ and $(f^\#E_{\tau-m})^*\in\Nev_{\pi\tau}^+$ simultaneously (where $f^\#$ varies in $\sgl{f,f^*}$).
By Krein's theorem \cite[p.38]{dB}, it is equivalent to say that $f^\# E_{\tau-m}$ is an entire function of exponential type at most $\pi\tau$. 
Furthermore, for $f\in \ex{-g}\Hil$ 
 $$\nm{f E_{\tau-m}}_2\simeq\nm{f\ex{\omega_{\tau-m}}}_2=\nm{f\ex{-\omega_m}}_2\simeq\nm{f(x)\ex{-\omega_\gamma(x+\eul)}}_2 \leq
  \nm{f\ex{-\omega_\gamma}}_2<\infty.$$  
In particular, $f^\# E_{\tau-m}$ belongs to the classical Paley--Wiener space $L^2_{\pi\tau}$, so the exponential type condition
 $|f^\#(z) E_{\tau-m}(z)|\leq C_\eps\ex{(\pi\tau+\eps)|z|}$ may be replaced with $|f^\#(z) E_{\tau-m}(z)|\leq C_\eps\ex{(\pi\tau+\eps)|y|}$ for $z=x+\eul y\in\cpl$.
Since $E_{\tau-m}$ is a Hermite--Biehler function, it is equivalent to say that this last inequation holds on $\sgl{\Im z\geq 0}$ only (where $f^\#$ varies in $\sgl{f,f^*}$).
There, $|E_{\tau-m}(z)|\simeq\ex{\pi\tau|y|}\ex{-\omega_m(z)}$, and hence the last inequation is equivalent to $|f^\#(z)|\ex{-\omega_m(z)}\leq C_\eps\ex{\eps |y|}$ for $\Im z\geq 0$,
 that is, to $|f(z)|\ex{-\omega_m(z)}\leq C_\eps\ex{\eps |y|}$ for all $z\in\cpl$.  
Therefore,
 $$\ex{-g}\Hil=\ens{f\mbox{ entire}}{\nm{f\ex{-\omega_\gamma}}_2<\infty,\ |f(z)|\ex{-\omega_m(z)}\leq C_\eps\ex{\eps |y|}}.$$
 
Let $\theta_0(x)=\omega_\gamma(x+\eul)-\omega_\gamma(x)$, so $\ex{-\omega_\gamma(x)}\simeq\ex{-\omega_m(x)}\ex{\theta_0(x)}$.
This extra weight is easily computed:
 $$\theta_0(x)=\int_{-\infty}^\infty\log\left|1-\frac{\eul}{t-x}\right|\gamma(t)\dif{t}
  =\dem\int_{-\infty}^\infty\log\left(1+\frac{1}{(t-x)^2}\right)\frac{\ph'(t)}{\pi}\dif{t},$$
where the integral is well-defined.
For a fixed $x\in\reel$, let $$A(t)=(\ph(t)-\ph(x))\log\left(1+\frac{1}{(t-x)^{2}}\right).$$ 
Notice that $\lim_{t\to x}A(t)=0$, while $\lim_{|t|\to\infty}A(t)=0$. 
Consequently, an integration by parts yields
 $$\theta_0(x)=\frac{1}{\pi}\int_{-\infty}^\infty\frac{1}{1+(t-x)^2}\,\frac{\ph(t)-\ph(x)}{t-x}\dif{t}.$$
The equation (\ref{BiLipschitz}) implies
 $$\theta_0(x)=\frac{1}{\pi}\int_{|t-x|<N}\frac{1}{1+(t-x)^2}\,\frac{\ph(t)-\ph(x)}{t-x}\dif{t} + O(1)$$
when $|x|\to\infty$. 
Since $\ph$ is non decreasing, it also implies that $|\ph(t)-\ph(x)|<2NC_2$ when $|t-x|<N$.
Therefore, 
 $$\theta_0(x)=\theta(x) - \frac{1}{\pi}\int_{|t-x|<N}\frac{t-x}{1+(t-x)^2}(\ph(t)-\ph(x))\dif{t} + O(1)
	            =\theta(x)+O(1),$$
where 
 $$\theta(x)=\frac{1}{\pi}\int_{|t-x|<N}\frac{\ph(t)-\ph(x)}{t-x}\dif{t}.$$
In conclusion, $\ex{\theta_0(x)}\simeq\ex{\theta(x)}$, and hence
 \begin{eqnarray*} \ex{-g}\Hil
  &=& \ens{f\mbox{ entire }}{\nm{f(x)\ex{-\omega_m(x)}\ex{\theta(x)}}_2<\infty,\ |f(z)|\ex{-\omega_m(z)}\leq C_\eps\ex{\eps|y|}}\\
  &=& \ens{f\in PW(m)}{\nm{f(x)\ex{-\omega_m(x)}\ex{\theta(x)}}_2<\infty},
 \end{eqnarray*}
with the above expression for $\theta(x)$.

\section{Examples: MC-spaces}\label{Sec:Examples}

In \cite{LS}, Lyubarskii and Seip studied structural properties of a larger class of de Branges spaces than the weighted Paley--Wiener spaces.
Let us define this larger class through a list of postulates \cite{P}.

A piecewise continuous function $\reel\to\reel$ is a \emph{mountain chain} if its graph consists of a succession of continuous pieces, called \emph{mountains} and \emph{plateaux}, satisfying the following conditions:
\begin{itemize}
 \item
  Each mountain has a Poissonian shape $\Kfrac{\eta}{(x-\xi)^2+\eta^2}$ with two sides and a summit, $(\xi,1/\eta)$;
 \item
  The bases of the mountains are $\simeq 1$;
 \item
  The summits have level more than 1; horizontally, they are bounded away from the endpoints of the mountain bases;
 \item
  The plateaux consist of horizontal segments of level $1$, without restriction on their lengths (finite or infinite).
\end{itemize}

Let $E$ be a Hermite--Biehler function without real zero and of phase $\ph$. 
For $x\in\reel$, let $\xi_x-\eul\eta_x$ be the zero of $E(z)$ closest to $x$ (with the smallest $x$-coordinate in case of equality). 
We postulate the existence of a $\delta>0$ such that the function
 \begin{equation*}
  \mu(x)=
   \left\{
    \begin{array}{ll}
     \Kfrac{\eta_x}{(x-\xi_x)^2+\eta_x^2}&\mbox{ if }\eta_x<\delta,\\
     1&\mbox{ otherwise}
    \end{array}
   \right.
 \end{equation*}
is a mountain chain, satisfying in addition
 $$\ph'(x)\simeq \mu(x).$$
Observe that each mountain of $\mu$ lies over a zero of $E(z)$ in the critical strip $-\delta < \Im z < 0$, and conversely.
We postulate that each such zero is simple. 
Finally, let $\sgl{(\xi_k,1/\eta_k)}_{k\in\rel}$ be an indexation of the summits of $\mu$ in order of $x$-coordinates, so that 
 $$\cdots<\xi_{-1}<\xi_0<\xi_1<\cdots$$
We postulate the following, weak limitation on the growth of the summits:
 $$|\log(\eta_k)-\log(\eta_l)|=O(|\xi_k-\xi_l|^{1-\eps})$$
uniformly in $l$ when $|k-l|\to \infty$, where $\eps>0$ is arbitrarily small.

In the sequel, a de Branges space $\Hil$ shall be called an \emph{MC-space} if $\Hil=\Hil(E)$ for a Hermite--Biehler function $E$ without real zero, whose phase
 satisfies the aformentioned, postulated properties.
In their profound study of weighted Paley--Wiener spaces, Lyubarskii and Seip proved the following theorem:

\begin{quote}\it
  If $\Hil$ is an MC-space, then there exists a real-entire $g(z)$ and a measurable, positive $m(x)\simeq 1$ such that $$\ex{g}PW(m)\subseteq\Hil.$$
  The majorants of these two spaces are both comparable with $\ex{g(x)}\ex{\omega_m(x)}$ on the real axis.
  In particular, if $\Hil$ is a weighted Paley--Wiener space, $\ex{g}PW(m)=\Hil$ (equality with equivalence of norms).
\end{quote}

\noindent For the full equality $\ex{g}PW(m)=\Hil(E)$ to hold, in addition to the properties of $E$ (more precisely, of $\ph'$) making $\Hil(E)$ an MC-space, the authors established the following, rather complicated criterion:

\begin{quote}\it
 For $\alpha\in\reel$, let $\Lambda_\alpha$ be the zero set of $\sin(\ph(x)-\alpha)$. 
 Let $g$ and $m$ be given by the previous theorem.
 Then, $\Hil(E)=\ex{g}PW(m)$ if and only if the following conditions are satisfied: Firstly, $\ex{-g(z)}\ex{-\omega_m(z)}E(z)$ is a (non entire) 
  function of exponential type on $\cpl$, which satisfies
   $$\int_{-\infty}^\infty\frac{\log^+|\ex{-g(x)}\ex{-\omega_m(x)}E(x)|}{x^2+1}\dif{x}<\infty.$$
 Secondly, for two $\alpha\in[0,\pi)$, $\Lambda_\alpha$ is separated, while 
  $$v(x)=\frac{\sin^2(\ph(x)-\alpha)}{\ph'(x)\dist{x}{\Lambda_\alpha}^2}$$
 satisfies the Muckenhoupt ($A_2$) condition.
\end{quote}

\noindent This last condition means that
 $$\Inv{|I|}\int_I v(x)\dif{x}\cdot\Inv{|I|}\int_I\Inv{v(x)}\dif{x}\lesssim 1$$
when $I$ ranges among the non empty finite intervals.

Notice that, if $\mu$ has arbitrarily high mountains (corresponding to zeroes of $E$ arbitrarily close to the 
 real axis), the relation $\ph'\simeq\mu$ ensures that $\ph$ is bi-Lipschitz for large distances.
However, $\ph$ fails to be bi-Lipschitz, since it contains arcs arbitrarily close to step functions (describing an almost
 horizontal segment, followed with an almost vertical one and an almost horizontal one, with angles arbitrarily close to the right angle).

The MC-spaces are best studied by comparing their Hermite--Biehler function $E$ with a simpler function, $F$, obtained from $E$ by 
 shifting down its zeroes in the critical strip $-\delta<\Im z<0$ to the axis $\Im z=-\delta$ \cite[p.991]{LS}.
$\Hil(E)$ is then related to $F$ as follows: for $\sigma(x)=\min(\eta_x,1)$, 
  $$\Hil(E)=\ens{f \mbox{ entire}}{\nm{f(x)\ph'(x)^\dem\sigma(x)^{-\dem}/F(x)}_2<\infty,\ f/F,\mbox{ and } f^*/F\in\Nev_0^+}.$$
In this particular context, $\ph'(x)^\dem\sigma(x)^{-\dem}$ is thus the analog of our extra-weight $\ex{\theta(x)}$.
Let us give a concrete example.

Let $\Hil(E)$ be the de Branges space generated by 
 $$E(z)=\,(z+\eul)\,\prod_{n=1}^\infty \left(1-\frac{z}{n-\eul n^{-\alpha}}\right)\left(1+\frac{z}{n+\eul n^{-\alpha}}\right),$$
where $0\leq \alpha\leq 2$.
Then, from the above factorization
 \begin{equation}\label{PhiPrime}
  \ph'(x)=\frac{1}{x^2+1}+\sum_{n\in\rel^*} \frac{|n|^{-\alpha}}{(x-n)^2+|n|^{-2\alpha}}.
 \end{equation}
Let $n_x$ be the integer closest to $x$ (less than $x$ in case of equality).
Then, 
 \begin{equation*}
  \mu(x)=
	 \left\{
	  \begin{array}{ll}
		 \kfrac{1}{(x^2+1)} &\mbox{ if } -\dem< x\leq\dem,\\
		 \kfrac{|n_x|^{-\alpha}}{(\,(x-n_x)^2+|n_x|^{-2\alpha})} & \mbox{ otherwise}
		\end{array}
	 \right.
 \end{equation*}
is a mountain chain.
In addition, $\ph'(x)-\mu(x)$ is positive and bounded above (when $|x|\to\infty$) by\\

\noindent $\displaystyle\frac{1}{x^2+1}+\!\!\!\sum_{|n-x|\geq\dem}\!\!\!\!\!'\ \ \frac{1}{|n|^\alpha(n-x)^2}
          \ \ \; \leq \ \frac{1}{x^2+1}+\!\!\!\!\!\sum_{|\,|n|-|x|\,|\geq\dem}\!\!\!\!\!\!\!\!\!'\ \ \ \frac{1}{|n|^\alpha(|n|-|x|)^2} + O(|x|^{-\alpha-2})$
\begin{eqnarray*} 
 &=& 2\;(\!\!\!\sum_{1\leq j<\frac{|x|}{2}}+\!\!\!\sum_{\frac{|x|}{2}\leq j\leq|x|-\dem}\!\!\!+\sum_{|x|+\dem\leq j}\!\!)\,\,\frac{1}{j^\alpha(j-|x|)^2} + O(|x|^{-2})\\
 &\lesssim& \!\frac{2^\alpha}{|x|^\alpha}\!\!\!\!\sum_{1\leq j<\frac{|x|}{2}}\frac{1}{j^\alpha(|x|-j)^{2-\alpha}} + 
            \left(\frac{2^\alpha}{|x|^\alpha}+\frac{1}{(|x|+\dem)^\alpha}\right)\sum_{k=0}^\infty\frac{1}{(\dem+k)^2}+ O(|x|^{-2})\\
 &\leq& \frac{2^\alpha}{|x|^\alpha}\sum_{1\leq j<\frac{|x|}{2}}\frac{1}{j^2} + O(|x|^{-\alpha}) 
 \ \; = \ \; O(|x|^{-\alpha}) \ \; \lesssim \ \; \mu(x).
\end{eqnarray*}

\noindent It follows that $\ph'(x)\simeq \mu(x)$, and hence, $\Hil(E)$ is an MC-space. 
By the above criterion, $\Hil(E)$ is indeed a weighted Paley--Wiener space for $0\leq \alpha<1$;
 if instead $1\leq\alpha<2$, the successive zeroes of $E$ approach the real axis too quickly for the Muckenhoupt~(A$_2$) condition to apply, so $\Hil(E)$ is not a
 weighted Paley--Wiener space (see \cite[p.1005]{LS} for a similar example).	
 
For $\Im z\gg 0$, the factorization of $E(z)$ and the simple relations 
 $$|z-n|<|z-(n-\eul|n|^{-\alpha})|<|z-(n-\eul)|$$
 ensure that $|\sin(\pi z)| \lesssim |E(z)| \lesssim |\sin(\pi (z+\eul))|.$
Therefore, $|E(z)|\simeq \ex{\pi \Im z} = \ex{\omega_1(z)}$ holds for $\Im z\gg 0$.
In this concrete case, our result thus applies for $m=1$ (without seeking another $m$ from the multiplier lemma).
Therefore, 
 $$\Hil(E)=\ens{f\in L^2_\pi}{\nm{f(x)\ex{\theta(x)}}_2<\infty}.$$ 
Let us study the extra-weight $\ex{\theta(x)}$ for $|x|>\dem$, with $N=\frac{5}{4}$.
Let
 $$\theta_0(x)=\frac{1}{\pi}\int_{-5/4}^{5/4}\frac{\tan^{-1}(t+x)-\tan^{-1}(x)}{t}\dif{t},$$
 $$\theta_n(x)=\frac{1}{\pi}\int_{-5/4}^{5/4}\frac{\tan^{-1}(|n|^\alpha(t+x-n))-\tan^{-1}(|n|^\alpha(x-n))}{t}\dif{t}.$$
The relation (\ref{PhiPrime}) then implies that $\theta(x)=\sum_{n\in\rel}\theta_n(x)$.
Observe that $\theta_0(x)$, $\theta_{n_x-1}(x)$, and $\theta_{n_x+1}(x)$ are bounded when $x$ varies in $\reel$. 
Let $\vartheta(x)=\theta_{n_x}(x)$.
It follows that $\theta(x)=\vartheta(x)+\sum_{|n-n_x|\geq 2}^{'}\theta_n(x)+O(1).$
By the mean value theorem, 
 $$\theta(x)-\vartheta(x) \leq \sum_{|n-x|\geq 3/2}\!\!\!\!\!\!\!'\,\,\,\,\,\,\,\frac{1}{\pi}\,\frac{|n|^{-\alpha}}{(|x-n|-\frac{5}{4})^2+|n|^{-2\alpha}}\int_{-5/4}^{5/4}\dif{t} + O(1) = O(1).$$
Thus, $\vartheta(x)$ stays at a bounded distance from $\theta(x)$, and hence $\ex{\vartheta(x)}\simeq\ex{\theta(x)}$.

Let us estimate $\vartheta(n+h)$ for $n\in\rel^*$ and $h\in[0,\dem]$ (the case where $h\in(-\dem,0)$ being similar). 
For $a<b$, let
 $$T(a,b)=\frac{1}{\pi}\int_a^b\frac{\tan^{-1}(\xi)-\tan^{-1}(|n|^\alpha h)}{\xi-|n|^\alpha h}\dif{\xi},$$
so $\vartheta(n+h)=T(-|n|^\alpha(\frac{5}{4}-h),|n|^\alpha(\frac{5}{4}+h))$. 
Notice that, when $n$ varies in $\rel^*$ and $h$ varies in $[0,\dem]$, both $T(-1,1)$ and $T(|n|^\alpha h - 1, |n|^\alpha h + 1)$ are bounded.
If $1\leq |n|^\alpha h - 1$, then $T(1,|n|^\alpha h-1)$ is also bounded: indeed, $\tan^{-1}(x)$ is concave downward on
 $[1,|n|^\alpha h]$, and hence 
 $$T(1,|n|^\alpha h-1)\ \leq \, \frac{1}{\pi}\,\frac{\tan^{-1}(|n|^\alpha h)-\tan^{-1}(1)}{|n|^\alpha h-1}
    \int_1^{|n|^\alpha h-1}\!\!\!\!\!\dif{\xi}\ \lesssim \, 1.$$
Therefore, 
 $$\vartheta(n+h)=T(-|n|^\alpha(\ff-h),-1)+T(|n|^\alpha h+1,|n|^\alpha(\ff+h))+O(1)$$
when $|n|\to\infty$, uniformly in $h\in[0,\dem]$.
The boundedness of $x(\frac{\pi}{2}-\tan^{-1}(x))$ when $x$ varies in $[0,\infty)$ then ensures that
 $$\frac{1}{\pi}\int_{-|n|^\alpha(\frac{5}{4}-h)}^{-1}\!\!\!\!\frac{\tan^{-1}(|n|^\alpha h)+\frac{\pi}{2}}{|n|^\alpha h-\xi}\,\dif{\xi} - T(-|n|^\alpha(\ff-h),-1) = O(1)$$
and 
 $$T(|n|^\alpha h+1,|n|^\alpha(\ff+h)) - \frac{1}{\pi}\int_{|n|^\alpha h+1}^{|n|^\alpha(\frac{5}{4}+h)}\frac{\frac{\pi}{2}-\tan^{-1}(|n|^\alpha h)}{\xi-|n|^\alpha h}\,\dif{\xi} = O(1)$$
when $|n|\to\infty$, uniformly in $h\in[0,\dem]$.
Therefore, 
 $$\vartheta(n+h)= \log(|n|^\alpha)-\frac{\tan^{-1}(|n|^\alpha h)+\frac{\pi}{2}}{\pi}\,\log(1+|n|^\alpha h) + O(1).$$
The boundedness of $x(\frac{\pi}{2}-\tan^{-1}(x))$ on $[0,\infty)$ finally ensures that 
 $$\log(1+|n|^\alpha h)-\frac{\tan^{-1}(|n|^\alpha h)+\frac{\pi}{2}}{\pi}\,\log(1+|n|^\alpha h)=O(1),$$
yielding in total
 $$\vartheta(n+h)=\log(|n|^\alpha)-\log(1+|n|^\alpha h)+O(1)$$
when $|n|\to\infty$, uniformly in $h\in [0,\dem]$.
A similar result holds for $h\in(-\dem,0)$.
In conclusion, for $x\notin(-\dem,\dem]$, 
 $$\ex{\theta(x)}\simeq\ex{\vartheta(x)}\simeq \frac{|n_x|^\alpha}{1+|n_x|^\alpha |x-n_x|}.$$
It is comparable to $\ph'(x)^{\dem}\sigma(x)^{-\dem}\simeq\mu(x)^{\dem}\sigma(x)^{-\dem}=|n_x|^\alpha/\sqrt{(x-n_x)^2|n_x|^{2\alpha}+1}$, as expected.

\paragraph{Concluding remark} 

We have seen that the MC-spaces share their structure with a larger class of spaces of entire functions, namely, the de Branges
 spaces with bi-Lipschitz phase for large distances. 
These last are of the form
 $$\Hil=\ens{f \mbox{ entire}}{\nm{f\ex{-\omega_m(x)}\ex{\theta(x)}}_2<\infty, \ |f(z)|\ex{-\omega_m(z)}\leq C_\eps\ex{\eps|\Im z|}}$$
for a measurable $m(x)\simeq 1$ and $\theta(x)=\int_{-N}^N \frac{\ph(t)-\ph(x)}{t-x}\dif{t}$, where $\ph(x)$ is a non decreasing, 
real-analytic function and $N>0$.
Contrary to the MC-spaces, they are defined by a simple condition on the phase of a corresponding Hermite--Biehler function, without
 assuming any knowledge of its zeroes.

Some de Branges spaces with bi-Lipschitz phase for large distances appear to be weighted Paley--Wiener spaces, despite an unbounded $\theta(x)$. 
The question is thus raised to find conditions on $\theta(x)$ for such a phenomenon to occur.
In this regards, it may be attempted to re-do the analysis of Lyubarskii and Seip by reasoning on the phase of $E$ (and its weight $\ex{\theta(x)}$) instead of the zeroes of $E$,
 and see if simplifications occur.

\paragraph{Acknowledgements}

The authors would like to acknowledge Prof.\@ Kristian Seip for reading this manuscript and giving encouraging comments; Prof.\@ Dmitry Jakobson for his kind invitation to present the result at \emph{McGill's Analysis Seminar}; and, Prof.\@ Marco Merkli for carefully reading this manuscript and improving its presentation.

\paragraph{Conflict of interests}

The authors declare that there is no conflict of interests regarding the publication of this paper.

\end{document}